\documentclass[12pt]{article}
\usepackage{amssymb}
\usepackage{amsfonts}
\usepackage{amsmath}
\usepackage{amsbsy}
\usepackage{latexsym}
\newtheorem {theorem}{Theorem}
\newtheorem {definition}{Definition}
\newtheorem {lemma}{Lemma}

\newtheorem {remark}{Remark}

\begin{document}

\title{Solutions of equations of viscous hydrodynamics via
stochastic perturbations of inviscid flows}

\author{Yuri E. Gliklikh\\
Mathematics Faculty, Voronezh State University\\
Universitetskaya pl. 1, Voronezh,
394006, Russia\\
yeg@math.vsu.ru}

\date{}

\maketitle


\begin{abstract}
\noindent We introduce a special stochastic perturbation of the flow
of diffuse matter as a curve in the group of diffeomorphisms of flat
$n$-dimensional torus such that the perturbed system yields a
solution of Burgers equation in the tangent space at unit of the
diffeomorphism group. The same perturbation of the flow of perfect
incompressible fluid yields a solution of Reynolds type equation but
under some special external force on the diffeomorphism group it
transforms into a solution of Navier-Stokes equation without
external force.
\end{abstract}

{\bf Keywords and phrases:} Group of diffeomorphisms; flat torus;
stochastic perturbation; diffuse matter; Burgers equation; perfect
incompressible fluid; Reynolds equation; Navier-Stokes equation.

{\bf 2000 Mathematics Subject Classification:} Primary 58J65, 60H30;
Secondary 76D05, 76D09

\section{Introduction \label{1}}
The paper is devoted to the Lagrangian approach to hydrodynamics in
terms of geometry of groups of diffeomorphisms, suggested for
perfect fluids by Arnold \cite{A} and Ebin and Marsden \cite{3}. In
previous works by the author it was found that the adequate
description of viscous fluids in this language requires involving
stochastic processes such that their expectations are flows of
viscous incompressible fluids (see, e.g., \cite{Gl1,Gl2}). In this
framework Newton's  second law on the groups of diffeomorphisms,
used in the case of perfect fluids, is replaced by its special
stochastic analogue in terms of Nelson's backward mean derivatives.
After transition to the tangent space at unit of diffeomorphisms
group, there arises the Navier-Stokes equation via a natural
modification of construction by Arnold, Ebin and Marsden that yields
the Euler equation in the case of perfect incompressible fluid. In
complete form this idea is realized in the case where the drift of
above-mentioned process on the group of diffeomorphisms is
right-invariant (see \cite{GG,GGG}).

Here we consider another case. We introduce a special stochastic
perturbation of a flow of diffuse matter, satisfying the
above-mentioned stochastic Newton's law, and show that the
corresponding curve in the tangent space at unit satisfies Burgers
equation. The same perturbation of a flow of perfect incompressible
flow without external force satisfies the above-mentioned stochastic
Newton's law as well, but yields a curve in the tangent space at
unit that is a solution of a Reynolds type equation. Nevertheless,
under the action of a certain special external force on the flow,
this curve becomes a solution of Navier-Stokes equation without
external force. We consider the fluid motion on the flat
$n$-dimensional torus (we recall the definition below).

The research was supported in part by RFBR Grants No. 07--01--00137
and 08--01--00155.

This paper originated from a discussion with B.~Rozovskii held long
ago at the University of Warwick. I am grateful to him for that
initial prompt.

\section{Preliminaries}\label{prelim}
Consider a stochastic process $\xi(t)$ in $\mathbb{R}^n $, where $t
\in [0,T]$, given on a certain probability space $(\Omega, {\cal F},
{\rm P})$ and such that $\xi(t)$ is an $L^1$-random variable for all
$t$. The ``present'' (``now'') for $\xi(t)$ is the least complete
$\sigma$-subalgebra ${\cal N}^\xi_t$ of $\cal F$ that includes
preimages of the Borel set of $\mathbb{R}^n $ under the map
$\xi(t):\Omega\rightarrow \mathbb{R}^n $. We denote by $E^{\xi}_t$
the conditional expectation with respect to ${\cal N}^\xi_t$. The
least complete $\sigma$-subalgebra that includes preimages of the
Borel set of $\mathbb{R}^n $ under all maps
$\xi(s):\Omega\rightarrow \mathbb{R}^n $ for $s\le t$ (resp. $s>t$)
is called the ``\textit{past}'' (resp. ``\textit{future}'')
$\sigma$-\textit{algebra} and is denoted by ${\cal P}^\xi_t$ (resp.
${\cal F}^\xi_t$).

Below we most often deal with the diffusion processes of the form
\begin{equation}
\xi (t) = \xi_0+\int_0^t  a(s,\xi (s))ds + \sigma w(t) \label{(8.6)}
\end{equation}
in $\mathbb{R}^n $ and flat torus ${\cal T}^n$ as well as natural
analogues of such processes on groups of diffeomorphisms
(infinite-dimensional manifolds). In (\ref{(8.6)}) $w(t)$ is a
\emph{Wiener process adapted to} $\xi(t)$, and $a(t,x)$ is a vector
field; $\sigma>0$ is a real constant.

Following Nelson (see, e.g., \cite{92} -- \cite{94}) we give the
next

\begin{definition}\label{8.1} (i) The \textit{forward mean
derivative} $D\xi(t)$ of the process $\xi(t)$ at $t$ is the
$L^1$-random variable of the form
\begin{equation}
D\xi (t)=\lim_{\Delta t \rightarrow +0}
E^\xi_t\left(\frac{\xi(t+\Delta t)-\xi(t)}{\Delta t}\right)
\label{(8.1)}
\end{equation}
where the limit is supposed to exist in $L^1(\Omega,{\cal F},{\rm
P})$ and $\Delta t \rightarrow +0$ means that $\Delta t\rightarrow
0$ and $\Delta t> 0$.

(ii) The \textit{backward mean derivative} $D_*\xi(t)$ of $\xi (t)$
at $t$ is the $L^1$-random variable
\begin{equation}
D_*\xi (t)=\lim_{\Delta t \rightarrow
+0}E^\xi_t\left(\frac{\xi(t)-\xi(t-\Delta t)}{\Delta t}\right)
\label{(8.2)}
\end{equation}
where (as well as in (i)) the limit is supposed to exist in
$L^1(\Omega,{\cal F},{\rm P})$ and $\Delta t \rightarrow + 0$ means
the same as in (i).\end{definition}

Notice that, generally speaking, $D\xi(t) \ne D_*\xi(t)$ (but, if
$\xi (t)$ almost surely (a.s.) has smooth sample trajectories, these
derivatives evidently coincide).

From the properties of conditional expectation it follows that $D\xi
(t)$ and $D_* \xi (t)$ can be represented as compositions of $\xi
(t)$ and Borel measurable vector fields
\begin{gather}
Y^0 (t,x) = \lim_{\Delta t \rightarrow +0} E\left(\frac{\xi
(t+\Delta t)-\xi
(t)}{\Delta t}\,\Big|\,{\xi (t)=x}\right),\notag\\
Y^0_* (t,x) = \lim_{\Delta t \rightarrow +0} E\left(\frac{\xi
(t)-\xi (t-\Delta t)}{\Delta t}\,\Big|\,{\xi (t)=x}\right)
\label{!!}
\end{gather}
on $\mathbb{R}^n $ (following \cite{Parthasarathy} we call them the
\emph{regressions}): $D\xi (t) = Y^0 (t, \xi (t))$ and $D_* \xi (t)
= Y^0_* (t, \xi (t))$.

\begin{lemma} \label{1.9} For a process of type $(\ref{(8.6)})$, we
have $D\xi(t)=a(t,\xi(t))$ and so $Y^0(t,x)=a(t,x)$\end{lemma}

For details of the proof, see, e.g., \cite{Gl1,Gl2}.

Let $Z (t, x)$ be $C^2$-smooth vector field on $\mathbb{R}^n $.

\begin{definition} \label{8.9} The $L^1$-limits of the form
\begin{equation}
DZ(t,\xi (t))=\lim_{\Delta t \rightarrow +0} E_t^\xi
\left({\frac{Z(t+\Delta t,\xi (t+\Delta t)) -Z(t,\xi (t))}{\Delta
t}}\right)  \label{(8.9)}
\end{equation}
and
\begin{equation}
D_*Z(t,\xi (t))=\lim_{\Delta t \rightarrow +0} E_t^\xi
\left({\frac{Z(t,\xi (t)) -Z(t-\Delta t,\xi (t-\Delta t))}{\Delta
t}}\right) \label{(8.10)}
\end{equation}
are called \textit{forward and backward}, respectively, \textit{mean
derivatives} of $Z$ along $\xi (\cdot)$ at time instant $t$.
\end{definition}

Certainly $DZ(t,\xi (t))$ and $D_* Z(t,\xi (t))$ can be represented
in terms of corresponding regressions, defined analogously to
(\ref{!!}). If it does not yield a confusion, we denote these
regressions by $DZ$ and $D_* Z$.

\begin{lemma} \label{8.10} For the process $(\ref{(8.6)})$ in
$\mathbb{R}^n$,
the following formulae take place:
\begin{gather}
DZ = {\frac{\partial}{\partial t}} Z+(Y^0\cdot \nabla)Z+
{\frac{\sigma^2}{2}}\nabla^2 Z,   \label{(8.11)}\\
\hspace{-145mm}\mbox{and}\notag\\
 D_*Z = {\frac{\partial}{\partial t}}
Z+(Y^0_* \cdot \nabla)Z- {\frac{\sigma^2}{2}}\nabla^2 Z
\label{(8.12)}
\end{gather}
where $\nabla=({\frac{\partial}{\partial
x^1}},...,{\frac{\partial}{\partial x^n}})$, $\nabla^2$ is the
Laplacian, the dot denotes the inner product in $\mathbb{R}^n $, the
vector fields $Y^0 (t, x)$ and $Y^0_* (t, x)$ are introduced in
(\ref{!!}).\end{lemma}

Recall that the $n$-dimensional torus ${\cal T}^n$ can be considered
as the quotient space of $\mathbb{R}^n$ with respect to the integral
lattice $\mathbb{Z}^n$. Introduce the Riemannian metric
$\langle\cdot,\cdot \rangle$ on ${\cal T}^n$ inherited form the
Euclidean inner product in $\mathbb{R}^n$. This metric is called
\textit{flat} and ${\cal T}^n$ with this metric is called the
\textit{flat torus}. Everywhere below we deal with the mean
derivatives and fluid motion on flat torus.

The main idea of description of viscous hydrodynamics in the
language of mean derivatives is as follows. Consider the vector
space $Vect^{(s)}$ of all Sobolev $H^s$-vector fields ($s>{\frac
n2}+1$) on ${\cal T}^n$.

Let a random flow $\xi (t,m)$ with initial data $\xi(0,m)=m\in {\cal
T}^n$ be given on a flat $n$-dimensional torus ${\cal T}^n$. Suppose
that it is a general solution of a stochastic differential equation
of the type
\begin{equation}
\xi (t,m)=m+ \int_0^t a(s,\xi (s,m))ds + \sigma w(t). \label{xi1}
\end{equation}
Let $D_*\xi(t,m)=v(t,\xi(t,m))$, where $v(t,m)$ is a $C^1$-smooth in
$t$ and $C^2$-smooth in $m\in {\cal T}^n$ vector field on ${\cal
T}^n$; let $\sigma >0$ be a real constant. Let $\xi (t,m)$ satisfy
the relation
\begin{equation}
D_* D_*\xi(t,m)=F(t,m), \label{Newtons4}
\end{equation}
where $F(t,m)$ is a vector field on ${\cal T}^n$. Taking into
account formula (\ref{(8.12)}), we obtain
\begin{gather}
D_* D_*\xi(t,m)=({\frac{\partial}{\partial
t}}v+(v,\nabla)v-{\frac {\sigma^2}2}\nabla^2 v)\notag 
={\frac{\partial}{\partial t}}v+(v,\nabla)v-{\frac
{\sigma^2}2}\nabla^2 v. \label{DDs4}
\end{gather}
Thus (\ref{Newtons4}) means that $v(t,m)$ satisfies the relation
\begin{equation}
{\frac{\partial}{\partial t}}v+(v,\nabla)v-{\frac
{\sigma^2}2}\nabla^2 v=F(t,m), \label{NSs4}
\end{equation}
that is the \emph{Burgers equation} with viscosity ${\frac
{\sigma^2}2}$ and external force $F(t,m)$. We interpret
\eqref{Newtons4} as a stochastic analogue of Newton's second law on
the group of Sobolev diffeomorphisms ${\cal D}^s({\cal T}^n)$ (see
the next Section).

The case of viscous incompressible fluids requires some additional
constructions. Introduce the $L^2$-inner product in $Vect^{(s)}$ by
the formula
\begin{equation}
(X,Y)=\int_{{\cal T}^n}\langle X(m),Y(m) \rangle\mu(dm) \label{L2}
\end{equation}
where $\langle \cdot,\cdot \rangle$ is the Riemannian metric on
${\cal T}^n$ and $\mu$ is the form of the Riemannian volume.

Denote by $\beta$ the subspace of $Vect^{(s)}$ consisting of all
divergence-free vector fields. Then consider the orthogonal
projection with respect to (\ref{L2}):
\begin{equation}
P: Vect^{(s)}\rightarrow \beta  \label{P}.
\end{equation}
 It follows from the Hodge decomposition that
the kernel of $P$ is the subspace consisting of all gradients. Thus,
for any $Y\in Vect^{(s)}$, we have
\begin{equation}
P(Y)=Y-{\rm grad} p, \label{P1}
\end{equation}
where $p$ is a certain $H^{s+1}$-function on ${\cal T}^n$, unique to
within the additive constant for a given $Y$.

Let a random flow $\xi (t,m)$ with initial data $\xi(0,m)=m\in {\cal
T}^n$ be given on a flat $n$-dimensional torus ${\cal T}^n$. Let
$\xi (t,m)$ be the general solution of a stochastic differential
equation of the type \eqref{xi1} and let
$D_*\xi(t,m)=u(t,\xi(t,m))$, where $u(t,m)$ is a $C^1$-smooth in $t$
and $C^2$-smooth in $m\in {\cal T}^n$ divergence-free vector field
on ${\cal T}^n$, and let $\sigma >0$ be a real constant. Suppose
that $\xi (t,m)$ satisfies the relation
\begin{equation}
P D_* D_*\xi(t,m)=F(t,m), \label{Newtons3}
\end{equation}
where $F(t,m)$ is a divergence-free vector field on ${\cal T}^n$.
Taking into account formulae (\ref{(8.12)}) and (\ref{P1}), we
obtain
\begin{gather}
P D_* D_*\xi(t,m)= P({\frac{\partial}{\partial
t}}u+(u,\nabla)u-{\frac {\sigma^2}2}\nabla^2 u)\notag\\
\qquad ={\frac{\partial}{\partial t}}u+(u,\nabla)u-{\frac
{\sigma^2}2}\nabla^2 u- {\rm grad} p. \label{DDs3}
\end{gather}
Thus (\ref{Newtons3}) means that $u(t,m)$ is divergence-free and
satisfies the relation
\begin{equation}
{\frac{\partial}{\partial t}}u+(u,\nabla)u-{\frac
{\sigma^2}2}\nabla^2 u- {\rm grad} p=F(t,m), \label{NSs3}
\end{equation}
that is the \emph{Navier-Stokes equation} with viscosity ${\frac
{\sigma^2}2}$ and external force $F(t,m)$.

We interpret (\ref{Newtons3}) as a stochastic analogue of Newton's
second law on the group of Sobolev diffeomorphisms ${\cal D}^s({\cal
T}^n)$ of the torus, subjected to the mechanical constraint.

\section{Basic notions from the geometry of groups of diffeomorphisms
of flat torus\label{2}}

The tangent bundle to ${\cal T}^n$ is trivial: $T{\cal T}^n={\cal
T}^n\times \mathbb{R}^n $. Note that the flat metric generates in
the second factor the inner product,  same as in the copy of
$\mathbb{R}^n$ from which the torus is obtained by factorization.

Consider the set ${\cal D}^s({\cal T}^n)$  of  all  diffeomorphisms
of ${\cal T}^n$ to itself belonging to the Sobolev space $H^s$,
where $s
> {\frac{1}{2}} n + 1$. Recall that for  $s > {\frac{1}{2}} n + 1$,
the maps belonging to $H^s$ class are $C^1$-smooth.

There  is  a  structure  of smooth (and separable) Hilbert manifold
on ${\cal D}^s({\cal T}^n)$ as  well  as the natural group
structures with respect to composition. A detailed description of
the structures and their interconnections can be found in \cite{3}.
Note that the tangent space $T_e{\cal D}^s({\cal T}^n)$ at the unit
$e=id$ is $Vect^{(s)}$ (see above). Recall that $T_e{\cal D}^s({\cal
T}^n)$ contains its subspace $\beta$ consisting of all
divergence-free vector fields on ${\cal T}^n$ belonging to $H^s$
(See Section \ref{prelim}).

The space $T_f{\cal D}^s({\cal T}^n)$, where $f \in {\cal D}^s({\cal
T}^n)$, consists  of  the maps $Y:{\cal T}^n \rightarrow TM$ such
that $\pi Y(m)=f(m)$, where $\pi:T{\cal T}^n \rightarrow {\cal T}^n$
is the natural projection. Obviously for any $Y\in T_f{\cal
D}^s({\cal T}^n)$, there exists unique $X\in T_e{\cal D}^s({\cal
T}^n)$ such that $Y=X\circ f$. In any $T_f{\cal D}^s({\cal T}^n)$ we
can define the $L^2$-inner product by analogy with (\ref{L2}) by the
formula
\begin{equation}
(X,Y)_f=\int_{{\cal T}^n}\langle X(m),Y(m) \rangle_{f(m)}\mu(dm).
\label{L22}
\end{equation}
The family of these inner products form the weak Riemannian metric
on ${\cal D}^s({\cal T}^n)$ (it generates the topology of the
functional space $H^0=L^2$, weaker than $H^s$).

The right translation $R_f: {\cal D}^s({\cal T}^n) \rightarrow {\cal
D}^s({\cal T}^n)$, where $R_f(\theta)= \theta\circ f$ for $\theta, f
\in {\cal D}^s({\cal T}^n)$, is $C^\infty$-smooth and thus  one may
consider right-invariant vector fields on ${\cal D}^s({\cal T}^n)$.
Note that the tangent to the right translation takes the form
$TR_fX=X\circ f$ for $X\in T{\cal D}^s({\cal T}^n)$.

The right-invariant vector field  $\bar X$ on ${\cal D}^s_\mu({\cal
T}^n)$ generated by a given vector $X \in T_e{\cal D}^s({\cal T}^n)$
is $C^k$-smooth if and only if the vector field $X$ on ${\cal T}^n$
is $H^{s+k}$-smooth. This fact is a consequence of the so-called
$\omega$-lemma (see \cite{3}) and it is valid also for more
complicated fields.  For example, if a tensor (or any other) field
on ${\cal T}^n$ is $C^\infty$-smooth, the corresponding
right-invariant field on ${\cal D}^s({\cal T}^n)$ is
$C^\infty$-smooth as well.

\begin{remark}\label{lx}
The left translation $L_f: {\cal D}^s({\cal T}^n) \rightarrow {\cal
D}^s({\cal T}^n)$, where $L_f(\theta)= \theta\circ f$ for $\theta, f
\in {\cal D}^s({\cal T}^n)$, is only continuous. Specify a vector
$x\in \mathbb{R}^n$ and denote by $l_x:{\cal T}^n\to{\cal T}^n$ the
diffeomorphism $l_x(m)=m+x$ {\em modulo factorization with respect
to the integral lattice}. Note that the left translation $L_{l_x}$
is $C^\infty$-smooth.
\end{remark}

Introduce the operators:
\begin{equation*}
B: T{\cal T}^n\rightarrow \mathbb{R}^n ,
\end{equation*}
the projection onto the second factor in ${\cal T}^n\times
\mathbb{R}^n $, and
\begin{equation}
A(m): \mathbb{R}^n \rightarrow T_m{\cal T}^n, \label{A}
\end{equation}
the converse to $B$ linear isomorphism of $\mathbb{R}^n $ onto the
tangent space to ${\cal T}^n$ at $m\in {\cal T}^n$. The map $A$ has
the following property. For the natural orthonormal frame  $b$ in
$\mathbb{R}^n $ we have an orthonormal frame  $A_m(b)$  in $T_m{\cal
T}^n$, the field of frames  $A(b)$ on ${\cal T}^n$ consists of
frames inherited from the constant frame  $b$. Thus, for a fixed
vector  $X \in \mathbb{R}^n $,  the vector field  $A(X)$  on ${\cal
T}^n$  is constant (i.e., it is obtained from the constant vector
field  $X$  on  $\mathbb{R}^n $ and has constant coordinates with
respect to $A(b))$ and, in particular,  $A(X)$  is $C^\infty$-smooth
and divergence-free since such is the constant vector field $X$ on
$\mathbb{R}^n $. So,  $A$  may be considered as a map  $A:
\mathbb{R}^n  \rightarrow \beta \subset T_e{\cal D}^s({\cal T}^n)$.

Introduce
\begin{equation}
Q_{g(m)}=A(g(m))\circ B, \label{(8)}
\end{equation}
where $g\in {\cal D}^s({\cal T}^n)$, $m\in {\cal T}^n$. For a vector
$Y\in T_f{\cal D}^s({\cal T}^n)$, we get $Q_gY=A(g(m))\circ
B(Y(m))\in T_g{\cal D}^s({\cal T}^n)$ for any $f\in {\cal D}^s({\cal
T}^n)$. In particular, $Q_eY\in Vect^{(s)}$. The operation $Q_e$ is
a formalization for ${\cal D}^s({\cal T}^n)$ of the usual
finite-dimensional operation that allows one to consider the
composition $X\circ f$ of a vector $X\in Vect^{(s)}$ and
diffeomorphism $f\in {\cal D}^s({\cal T}^n)$ as a vector in
$Vect^{(s)}$. It denotes the shift of a vector, applied at the point
$f(x)$, to the point $x$ with respect to global parallelism of the
tangent bundle to torus.

\begin{lemma}[\cite{G3}]\label{L1}
The following relations hold:
\begin{gather}
TR_{g^{-1}}(Q_gX)=Q_e(TR_{g^{-1}}X); \label{**i}\\
\hspace{-145mm}\mbox{and}\notag\\
 TR_g(Q_{g^{-1}}X)= Q_e(TR_gX).
\label{**ii}
\end{gather}
\end{lemma}

{\bf Proof.}  By the above formulae we see that $Q_e(TR_{g^{-1}}X)$
sends the point $m\in {\cal T}^n$ to $(m,X(g^{-1}(m)))$. On the
other hand, $Q_gX=(g(m), X(m))$, and hence
\[
 TR_{g^{-1}}(Q_gX)=(m,X(g^{-1}(m)))
\]
so that (\ref{**i}) is proved. Formula (\ref{**ii}) follows from
(\ref{**i}) under the replacement of $g$ by $g^{-1}$. $\Box$

It should be pointed out that $Q_g$ is the global parallel
translation at $g$ on ${\cal D}^s({\cal T}^n)$ generated by global
parallelism on ${\cal T}^n$. It turns out that $Q_g$ is the
parallelism of Levi-Civita connection of metric \eqref{L22}. Thus,
for a smooth vector field $Y(t)$ along a smooth curve $g(t)$ in
${\cal D}^s({\cal T}^n)$, the covariant derivative ${\frac{\bar
D}{dt}}Y(t)$ at a time instant $t^*$ is defined as
\begin{equation}\label{barDdt}
{\frac{\bar
D}{dt}}Y(t)_{|t=t^*}={\frac{d}{dt}}(Q_{g(t^*)}Y(t))_{|t=t^*}.
\end{equation}

As usual, a smooth curve $g(t)$ in ${\cal D}^s({\cal T}^n)$ such
that
\begin{equation}\label{NewtonDs}
{\frac{\bar D}{dt}}\dot g(t)=0,
\end{equation}
is called \emph{geodesic}. In the framework of Lagrangian approach
to hydrodynamics \cite {A,3}, $g(t)$ describes the motion of
so-called {\em diffuse matter} on ${\cal T}^n$ without external
force (the case with non-zero external force is constructed in
analogy with that of perfect incompressible fluid below). For such
$g(t)$ introduce the vector $v(t)\in T_e{\cal D}^s({\cal T}^n)$
(i.e., the $H^s$-vector field $v(t,m)$ on ${\cal T}^n$) by the
formula $v(t)=\dot g(t)\circ g^{-1}(t)$. It is shown that $v(t)$
satisfies the \emph{Hopf equation} (sometimes called \emph{Burgers
equation without viscous term}):
\begin{equation}\label{Hopf}
{\frac{\partial}{\partial t}}v+(v\cdot\nabla)v=0.
\end{equation}

\begin{remark}\label{rightinv}
It is shown in \cite{3} that if $g(t)$ is a geodesic, then for every
$f\in {\cal D}^s({\cal T}^n)$, the curve $R_fg(t)$ is also a
geodesic.
\end{remark}

\begin{lemma}\label{Qlx}
Let $g(t)$ satisfy \eqref{NewtonDs} and $x\in\mathbb{R}^n$ be an
arbitrary specified vector. Then $l_xg(t)$ satisfies
\eqref{NewtonDs} as well, where $l_x$ is introduced in Remark
\ref{lx}.
\end{lemma}

{\bf Proof.}  Note that the parallel translations given by operators
$Q_g$ and $Q_f$ (see \eqref{(8)}), are commutative. By construction
${\frac{d}{dt}}l_xg(t)={\frac{d}{dt}}(g(t,m)+x)=Q_{l_x}\dot g(t)$
and $Q_{l_x}(Q_{g(t^*)}\dot g(t))=Q_{l_xg(t^*)}{\frac{d}{dt}}l_x
g(t)$. Taking into account definition \eqref{barDdt} of covariant
derivative we obtain ${\frac{\bar D}{dt}}{\frac {d}{dt}}l_x
g(t)=Q_{l_x}{\frac{\bar D}{dt}}\dot g(t)=0.$ $\Box$

Introduce the subspace $\beta_f\subset T_f{\cal D}^s({\cal T}^n)$ as
$TR_f\beta$ with $\beta$ introduced in \S \ref{prelim} (see, e.g.,
\eqref{P}). Having done this at every $f\in {\cal D}^s({\cal T}^n)$,
we obtain the smooth subbundle $\bar \beta$ of $T{\cal D}^s({\cal
T}^n)$ that in constructions below will be considered as constraint.
This constraint is holonomic, i.e., the distribution $\bar\beta$ is
integrable. The integral manifold going through $e$ is the
submanifold and subgroup ${\cal D}^s_\mu({\cal T}^n)$ in ${\cal
D}^s({\cal T}^n)$ that consists of $H^s$-diffeomorphisms preserving
the volume (see details in \cite{3}).

Notice that for $Y\in \beta_f$ the vector $Q_eY$ may not belong to
$\beta_e=\beta$.

Consider the map ${\bar P}: T{\cal D}^s({\cal T}^n)\rightarrow
\bar\beta$ determined for each  $f \in {\cal D}^s({\cal T}^n)$ by
the formula
\begin{equation*}
{\bar P}_f = TR_f\circ P\circ TR_f^{-1},
\end{equation*}
where $P=P_e:Vect^{(s)}=T_e{\cal D}^s({\cal T}^n)\rightarrow
\beta=\beta_e=T_e{\cal D}^s_\mu({\cal T}^n)$ is the orthogonal
projection introduced in (\ref{P}). It is obvious that $\bar P$ is
${\cal D}_\mu^s({\cal T}^n)$-right-invariant. An important and
rather complicated result (see \cite{3}) states that $\bar P$ is
$C^\infty$-smooth.

It is a routine fact of differential geometry that the covariant
derivative ${\frac{\widetilde D}{dt}}Y(t)$ of a vector field $Y(t)$
along a curve $g(t)$ in ${\cal D}^s_\mu({\cal T}^n)$ is defined by
the relation \begin{equation*} {\frac{\widetilde D}{dt}}Y(t)=\bar
P{\frac{\bar D}{dt}}Y(t).
\end{equation*}

Let $\bar F(t,g,Y)$, $Y\in T_g{\cal D}^s_\mu({\cal T}^n)$, be a
(force) vector field on ${\cal D}^s_\mu({\cal T}^n)$. Consider a
curve $g(t)$ satisfying the equation
\begin{equation}\label{Euler1}
{\frac{\widetilde D}{dt}}\dot g(t)=\bar F(t,g(t),\dot g(t)).
\end{equation}
Denote by $u(t)$ the curve in $T_e{\cal D}^s_\mu({\cal T}^n)$ (i.e.,
a divergence free vector field on ${\cal T}^n$) obtained by right
translations of vectors $\dot g(t)$, i.e., $u(t)=\dot g(t)\circ
g^{-1}(t)=TR_{g(t)}^{-1}\dot g(t)$. It is shown in \cite{3} that
$u(t)$ satisfies the Euler equation
\begin{equation}\label{Euler2}
{\frac{\partial}{\partial t}}u+(u\cdot\nabla)u-{\rm grad}
p=TR_g^{-1} \bar F(t,g(t),u(t,g(t))).
\end{equation}

It should be pointed out that \eqref{Euler1} is Newton's second law
with force $\bar F$ that describes the motion of perfect
incompressible fluid on ${\cal T}^n$ under the action of force
$$TR_g^{-1} \bar F(t,g(t),u(t,g(t)))$$ depending of the
``configuration of fluid''. Recall that a curve satisfying
\eqref{Euler1} with $\bar F=0$ is a geodesic.

\begin{remark}\label{rightinv2}
It is shown in \cite{3} that if $g(t)$ is a geodesic in ${\cal
D}^s_\mu({\cal T}^n)$, then for every $f\in {\cal D}^s_\mu({\cal
T}^n)$ the curve $R_fg(t)$ is also a geodesic, but nothing like
Lemma \ref{Qlx} is valid in this case.
\end{remark}

If $\bar F$ is a right-invariant vector field on ${\cal
D}^s_\mu({\cal T}^n)$ such that $\bar F_e=F$, where $F$ is a
divergence free vector field on ${\cal T}^n$, then \eqref{Euler2}
turns into
\begin{equation}\label{Euler3}
{\frac{\partial}{\partial t}}u+(u\cdot\nabla)u-{\rm grad} p=F.
\end{equation}

Consider the map $\bar A: {\cal D}^s({\cal T}^n) \times \mathbb{R}^n
\rightarrow T{\cal D}^s({\cal T}^n)$ such that $\bar A_e:
\mathbb{R}^n  \rightarrow \beta_e=T_e{\cal D}^s_\mu({\cal
T}^n)\subset T_e{\cal D}^s({\cal T}^n)$  is equal to $A$ introduced
by \eqref{A}, and for every $g \in {\cal D}^s({\cal T}^n)$, the map
$\bar A_g: \mathbb{R}^n \rightarrow T_g{\cal D}^s({\cal T}^n)$  is
obtained from $\bar A_e$ by means of the right translation, i.e.,
for $X\in \mathbb{R}^n$:
\begin{equation}
\bar A_g(X) = T R_g\circ A_e(X) = (A\circ g)(X). \label{barA}
\end{equation}
Since  $A$  is  $C^\infty$-smooth, it follows from $\omega$-lemma
that $\bar A$ is $C^\infty$-smooth jointly in  $X \in \mathbb{R}^n $
and $g \in {\cal D}^s({\cal T}^n)$. In particular, the restriction
$\bar A: {\cal D}^s_\mu({\cal T}^n) \times \mathbb{R}^n  \rightarrow
T{\cal D}^s_\mu({\cal T}^n)$ is $C^\infty$-smooth and the
right-invariant vector field $\bar A(X)$ is $C^\infty$-smooth on
${\cal D}_\mu^s({\cal T}^n)$ for every specified $X \in \mathbb{R}^n
$.

For any point $m\in {\cal T}^n$, denote by $exp_m:T_m{\cal
T}^n\rightarrow {\cal T}^n$ the map that sends the vector $X\in
T_m{\cal T}^n$ into the point $m+X$ in ${\cal T}^n$, where $m+X$ is
obtained {\em modulo factorization with respect to the integral
lattice}, i.e., by the following procedure: We take a certain point
in $\mathbb{R}^n$ corresponding to $m\in {\cal T}^n$, (denote it
also by $m$) and $X\in \mathbb{R}^n=T_m\mathbb{R}^n$, then we
identify $\mathbb{R}^n$ with $T_m\mathbb{R}^n=\mathbb{R}^n$, find
$m+X$ in $\mathbb{R}^n$ and pass from $\mathbb{R}^n$ to ${\cal T}^n$
by factorization with respect to $\mathbb{Z}^n$. The field of maps
$exp$ at all points generates the map $\overline{exp}:T_e{\cal
D}^s({\cal T}^n)\rightarrow {\cal D}^s({\cal T}^n)$ that sends the
vector $X\in T_e{\cal D}^s({\cal T}^n)$ (i.e., a vector field on
${\cal T}^n$) to $e+X\in {\cal D}^s({\cal T}^n)$, where $e+X$ is the
diffeomorphism of ${\cal T}^n$ of the form $(e+X)(m)=m+X(m)$.

Consider the composition $\overline{exp}\circ\bar
A_e:\mathbb{R}^n\rightarrow {\cal D}^s({\cal T}^n)$. By the
construction of $\bar A_e$ for any $X\in \mathbb{R}^n$ we get
$\overline{exp}\circ\bar A_e(X)(m)=m+X$, i.e., the same vector $X$
is added to every point $m$. Thus, obviously,
$\overline{exp}\circ\bar A_e(X)\in {\cal D}^s_\mu({\cal T}^n)$ and
so $\overline{exp}\circ\bar A_e$ sends $\mathbb{R}^n$ to ${\cal
D}^s_\mu({\cal T}^n)$.

Let $w(t)$ be a Wiener process in $\mathbb{R}^n$ defined on a
certain probability space $(\Omega, {\cal F}, {\sf P})$. Introduce
the process
\begin{equation}
W^{(\sigma)}(t)=\overline{exp}\circ\bar A_e(\sigma w(t)) \label{W}
\end{equation}
in ${\cal D}^s_\mu({\cal T}^n)$, where $\sigma > 0$ is a real
constant. By construction, for $\omega\in\Omega$ the corresponding
sample trajectory $W^{(\sigma)}_\omega(t)$ is the diffeomorphism of
the form $W^{(\sigma)}_\omega(t)(m)=m+\sigma w_\omega(t)$ so that
the same sample trajectory $\sigma w_\omega(t)$ of $w(t)$ is added
to each point $m\in {\cal T}^n$. In particular, this clarifies the
fact that $W^{(\sigma)}(t)$ takes values in ${\cal D}^s_\mu({\cal
T}^n)$.

\begin{remark}\label{ttt}
Note that for a ``specified $\omega\in \Omega$'' (i.e., a.s. for
$\omega\in\Omega$) and for specified $t\in\mathbb{R}$ the value of
$\sigma w_\omega(t)$ is a constant vector in $\mathbb{R}^n$. Then,
for ``given'' $\omega$ and $t$, the action of
$W^{(\sigma)}_\omega(t)$ coincides with that of $l_{\sigma
w_\omega(t)}$ (see Remark \ref{lx}).
\end{remark}

\section{Viscous hydrodynamics}

In this section we take $s>{\frac{n}{2}}+2$ so that the
diffeomorphisms from ${\cal D}^s({\cal T}^n)$ and ${\cal
D}^s_\mu({\cal T}^n)$ are $C^2$-smooth and $Vect^{(s)}$ consists of
$C^2$-smooth vector fields.

Everywhere below we use the same process $W^{(\sigma)}$ constructed
from a specified Wiener process $w(t)$ in $\mathbb{R}^n$ by formula
\eqref{W}. If, in the formula, there are several random elements
with subscript $\omega$, this means that they all are taken at ``the
same'' $\omega\in\Omega$, i.e., sometimes the formula may be
considered as description of a non-random element depending on the
parameter $\omega\in\Omega$.

Let $g(t)$ be a solution of \eqref{NewtonDs} with initial conditions
$g(0)=e$ and $\dot g(0)=v_0\in T_e{\cal D}^s({\cal T}^n)$. It is
shown in \cite{3} that such a solution exists on a certain time
interval $t\in [0,T]$ (for the sake of convenience we take a closed
interval inside the domain of $g(t)$). Recall that $g(t)$ is a flow
of diffuse matter without external forces. Consider $v(t)=\dot
g(t)\circ g^{-1}(t)\in T_e{\cal D}^s({\cal T}^n)$. This
infinite-dimensional vector considered as a vector field on ${\cal
T}^n$, will be also denoted by $v(t,m)$. Recall that this vector
field satisfies the Hopf equation \eqref{Hopf}.

Consider a process on ${\cal D}^s({\cal T}^n)$ of the form
$\eta(t)=W^{(\sigma)}(t)\circ g(t)$, $t\in [0,T]$, where
$W^{(\sigma)}(t)$ is introduced in \eqref{W}. In finite-dimensional
notation, $\eta(t)$ is a random diffeomorphism of ${\cal T}^n$ of
the form $\eta(t,m)=g(t,m)+\sigma w(t)$ modulo the factorization
with respect to integral lattice. Introduce the process
$\xi(t)=\eta(T-t)$, i.e., in finite-dimensional notation,
$\xi(t,m)=g(T-t,m)+\sigma w(T-t)$.

Since $w(t)$ is a martingale with respect to its own ``past'', one
can easily derive from the properties of conditional expectation
that $D_*\xi(t)=\dot g(T-t,m)=v(T-t, g(T-t,m))$ and so $\bar
PD_*D_*\xi(t)={\frac{\bar D}{ds}}\dot g(s)_{|s=T-t}=0$.

Consider the random process \[\xi_t(s)=\xi(s)\circ
\xi^{-1}(t)=W^{(\sigma)}(T-s)\circ g(T-s)\circ g^{-1}(T-t)\circ
(W^{(\sigma)}(T-t))^{-1}.\] Notice that the random diffeomorphism
$(W^{(\sigma)}(t))^{-1}$ acts by the rule
\[(W^{(\sigma)}(t))^{-1}(m)=m-\sigma w(t).\] Obviously $\xi_t(t)=e$.
The finite-dimensional description of this process can be given as
follows.

By construction \[m=\xi(t,\xi^{-1}(t,m))=g(T-t,\xi^{-1}(t,m))+\sigma
w(T-t).\] Then \[g(T-t,\xi^{-1}(t,m))=m-\sigma w(T-t)\] and so
\[\xi^{-1}(t,m)=g^{-1}(T-t, m-\sigma w(T-t)).\] Thus,
\[\xi_t(s,m)=\xi(s,g^{-1}(T-t,m-\sigma w(T-t))=\]\[g\bigl(T-s,g^{-1}(T-t,
m-\sigma w(T-t))\bigr)+\sigma w(T-s).\] Notice that
$\xi_t(t,m)=m-\sigma w(T-t)+\sigma w(T-t)=m$, i.e., indeed,
$\xi_t(t)=e$ on ${\cal D}^s({\cal T}^n)$.

Since $\xi_t(t)=e$, the ``present'' $\sigma$-algebra ${\cal
N}_t^{\xi_t}$ is trivial and so the conditional expectation with
respect to it coincides with ordinary mathematical expectation.
Hence, using the relation between $v(t)$ and $g(t)$ and definition
of $D_*$, one can easily derive that
\begin{gather}\label{HopfBurgers}
D_*\xi_t(s)_{|s=t}=E\bigl(v(T-t,m-\sigma w(T-t))\bigr)=
E\bigl(Q_eTR_{W^{(\sigma)}(T-t)}^{-1} v(T-t)\bigr)
\end{gather}
(note that here $t$ is specified and the derivative is taken with
respect to $s$).

Introduce on ${\cal T}^n$ the vector field $V(t,m)=
E\bigl(v(t,m-\sigma w(t))\bigr)$. We also denote this field, as an
infinite dimensional vector, by
$V(t)=E\bigl(Q_eTR_{W^{(\sigma)}(t)}^{-1} v(t)\bigr)$. Formula
\eqref{HopfBurgers} means that $D_*\xi_t(s)_{|s=t}=V(T-t)$.

\begin{theorem}\label{Burgers1}
The vector field $V(T-t,m)$ satisfies the Burgers equation
\begin{equation}\label{Burgers2}
{\frac{d}{dt}}V(T-t,m)+(V(T-t,m)\cdot\nabla)V(T-t,m)-
{\frac{\sigma^2}{2}}\nabla^2V(T-t,m)=0,
\end{equation}
where $\nabla^2$ is the Laplace-Beltrami operator which on the flat
torus coincides with the ordinary Laplacian.
\end{theorem}

{\bf Proof.} For $t\in[0,T]$ and $\omega\in\Omega$, introduce the
curve $\zeta_{t,\omega}(s)$ in $s\in[0,T]$ depending on parameter
$\omega$, by the formula
\[
\zeta_{t,\omega}(s)=R^{-1}_{W^{(\sigma)}_\omega(T-t)}g(T-s,g^{-1}(T-t))=g(T-s,g^{-1}(T-t,m-\sigma
w(T))).
\]
Note the only difference between $\zeta_{t,\omega}(s)$ and
$\xi_{t,\omega}(s)$: there is the stochastic summand $\sigma w(T-s)$
in the expression for $\xi_{t,\omega}(s)$ while it is absent in
$\zeta_{t,\omega}(s)$. This means, in particular, that
$\zeta_{t,\omega}(s)$ is an a.s. smooth curve with random initial
condition \[\zeta_{t,\omega}(t)=(W^{(\sigma)}_\omega(T-t))^{-1}.\]

Note that
${\frac{d}{ds}}\zeta_{t,\omega}(s)_{|s=t}=-TR^{-1}_{W^{(\sigma)}_\omega(T-t)}v(T-s)$.
Since $g(T-s)$ is a geodesic, from Remark \ref{rightinv} it follows
that ``for almost all specified $\omega$'' (i.e., a.s. for $\omega
\in\Omega$) the curve $\zeta_{t,\omega}(s)$ is also a geodesic,
i.e., ${\frac{\bar D}{ds}}{\frac{d}{ds}} \zeta_{t,\omega}(s)=0$. By
Remark \ref{ttt} the action of the diffeomorphism
$W^{(\sigma)}_\omega(t)$ coincides with that of $l_{\sigma
w_\omega(t)}$. Hence, by Lemma \ref{Qlx} the curve
\[(W^{(\sigma)}_\omega(T-t))\zeta_{t,\omega}(s)=l_{(\sigma
w_\omega(T-t))}\zeta_{t,\omega}(s)\] is a.s. geodesic as well, i.e.,
${\frac{\bar D}{ds}}{\frac{d}{ds}}l_{(\sigma w_\omega(T-t))}
\zeta_{t,\omega}(s)=0$. Note that
\[{\frac{d}{ds}}l_{(\sigma w_\omega(T-t))}
\zeta_{t,\omega}(s)_{|s=t}=Q_e{\frac{d}{ds}}
\zeta_{t,\omega}(s)_{|s=t}=-Q_eTR_{W^{(\sigma)}_\omega(T-t)}^{-1}
v(T-t).\] Recall that $EQ_eTR_{W^{(\sigma)}_\omega(T-t)}^{-1}
v(T-t)=V(T-t)$ and $D_*\xi_t(s)_{|s=t}=V(T-t)$ (see above). Then
from the above arguments and construction we derive that
\[
D_*D_*\xi_t(s)_{|s=t}=D_*V(T-t,\xi_t(s))_{|s=t}=-E\left({\frac{\bar
D}{ds}}{\frac{d}{ds}}l_{(\sigma
w_\omega(t))}\zeta_{t,\omega}(s)_{|s=t}\right)=0.
\]
But since $D_*\xi_t(s)_{|s=t}=V(T-t)$,  by formula \eqref{(8.12)}
the backward derivative $D_*V(T-t,\xi_t(s))$ coincides with the
left-hand side of \eqref{Burgers2}. Hence \eqref{Burgers2} is
fulfilled. $\Box$

Now, let us turn to the case of incompressible fluids. Let $g(t)$ be
a solution of \eqref{Euler1} on ${\cal D}^s_\mu({\cal T}^n)$ with
$\bar F=0$ and with initial conditions $g(0)=e$ and $\dot
g(0)=u_0\in T_e{\cal D}^s_\mu({\cal T}^n)$. As well as in the case
of diffuse matter, it is shown in \cite{3} that such a solution
exists in a certain time interval $t\in [0,T]$ (for the sake of
convenience we again take a closed interval inside the domain of
$g(t)$). Recall that $g(t)$ is a flow of perfect incompressible
fluid without external forces. Consider $u(t)=\dot g(t)\circ
g^{-1}(t)\in T_e{\cal D}^s_\mu({\cal T}^n)$. This
infinite-dimensional vector considered as a divergence free vector
field on ${\cal T}^n$, will be denoted $u(t,m)$. Recall that this
vector field satisfies the Euler equation \eqref{Euler3} without
external forces (see \S \ref{2}).

Since $W^{(\sigma)}(t)$ takes values in ${\cal D}^s_\mu({\cal T}^n)$
(see \S \ref{2}), we can repeat on ${\cal D}^s_\mu({\cal T}^n)$ the
above constructions for ${\cal D}^s({\cal T}^n)$, i.e., introduce
$\eta(t)=W^{(\sigma)}(t)\circ g(t)$, where $t\in [0,T]$, and
$\xi(t)=\eta(T-t)$ (i.e., in finite-dimensional notation
$\xi(t,m)=g(T-t,m)+\sigma w(T-t)$). One can easily see that
$D_*\xi(t)=\dot g(T-t,m)=u(T-t, g(T-t,m))$ and so ${\bar
PD_*D_*\xi(t)={\frac{\widetilde D}{ds}}\dot g(s)_{|s=T-t}=0}$ on
${\cal D}^s_\mu({\cal T}^n)$.

As well as above the process $\xi_t(s)=\xi(s)\circ \xi^{-1}(t)$ has
the property $\xi_t(t)=e$. Its finite-dimensional description is
quite analogous to the case of ${\cal D}^s({\cal T}^n)$.

Introduce on ${\cal T}^n$ the vector field $U(t,m)=
E\bigl(u(t,m-\sigma w(t))\bigr)$ (a direct analog of $V(t,m)$). We
also denote this field as an infinite dimensional vector by
$U(t)=E\bigl(Q_eTR_{W^{(\sigma)}(t)}^{-1} u(t)\bigr)$.

\begin{lemma}\label{Udivfree}
The vector field $U(t,m)$ is divergence free.
\end{lemma}

{\bf Proof.} By construction, for an elementary event
$\omega\in\Omega$, $(W^{(\sigma)}(t)_\omega)^{-1}$ is a shift of the
entire torus by a constant vector. Hence, $Q_e$ applied to
$TR_{W^{(\sigma)}(t)_\omega}^{-1}u(t)$ means the parallel
translation on torus of the entire divergence free vector field
$u(t)$ by the same constant vector back. Thus
$Q_eTR_{W^{(\sigma)}(t)}^{-1}u(t)$ is a random divergence free
vector field on the torus. Hence its expectation is divergence free.
$\Box$

So, $U(t)\in T_e{\cal D}^s_\mu({\cal T}^n)$. In particular, we have
proved above that
\begin{equation}\label{D*=U}
D_*\xi_t(s)_{|s=t}=U(T-t).
\end{equation}
Since nothing like Lemma \ref{Qlx} holds on ${\cal D}^s_\mu({\cal
T}^n)$, we have
$PD_*D_*\xi_t(s)_{|s=t}=D_*U(T-t,\xi_t(s))_{|s=t}\neq 0$ (the value
of this mean derivative is calculated in Remark \ref{4.1} below).
Hence there is no analogue of Theorem \ref{Burgers1}. We can prove
only the following

\begin{theorem}\label{mainaux}
Vector field $U(t,m)$ satisfies the following Reynolds type equation
(see, e.g., $\cite{Sedov}$)
\begin{equation}\label{Reynolds}
{\frac{\partial}{\partial t}}U+
E\Bigl[\Bigl((u\cdot\nabla)u\Bigr)(t,m-\sigma w(t))\Bigr]
-{\frac{\sigma^2}{2}}\nabla^2 U-{\rm grad}\,p=0.
\end{equation}
\end{theorem}

{\bf Proof.}  It follows from It\^o formula that
\begin{gather*}
du(t,m-\sigma w(t))={\frac{\partial u}{\partial t}}(t,m-\sigma
w(t))dt +{\frac{\sigma^2}{2}}\nabla^2u(t,m-\sigma w(t))dt
-\sigma\,u'\, dw(t),
\end{gather*}
where $\nabla^2$, as well as above, is the Laplace-Beltrami operator
and $u'$ is the linear operator of derivative of $u$ in $m\in {\cal
T}^n$.

Recall that $u(t,m)$ satisfies the Euler equation without external
force, i.e., ${\frac{\partial u}{\partial t}}=-P((u\cdot\nabla)u)$.
Since
\[E({\frac{d}{d t}}u(t,m-\sigma w(t)))={\frac{\partial}{\partial
t}}Eu(t,m-\sigma w(t))={\frac{\partial}{\partial t}}U(t)\] and
$E\bigl(\sigma(\nabla u) dw(t)\bigr)=0$, we derive that
\begin{gather*}
{\frac{\partial}{\partial t}}U=E\left({\frac{d}{dt}}u(t,m-\sigma
w(t))\right)=\\E\Bigl[-P\Bigl((u\cdot\nabla)u\Bigr)(t,m-\sigma
w(t))+{\frac{\sigma^2}{2}}\nabla^2 u(t,m-\sigma
w(t))\Bigr]=\\-E\Bigl[\Bigl((u\cdot\nabla)u\Bigr)(t,m-\sigma
w(t))\Bigr] +{\frac{\sigma^2}{2}}\nabla^2 U+{\rm grad}\,p.
\end{gather*}
So, \eqref{Reynolds} is satisfied. $\Box$

There are usual methods for transforming \eqref{Reynolds} into the
standard Reynolds form (see \cite{Sedov}). For a divergence-free
vector field $X(m)$ on ${\cal T}^n$ (i.e., for a vector $X\in
T_e{\cal D}^s_\mu({\cal T}^n)$) introduce the random divergence free
vector field
\[\breve U_X(t,m)=X(m-\sigma w(t))-E(X(m-\sigma w(t))\] (i.e.,
the vector \[\breve U_X(t)=Q_eTR_{W^{(\sigma)}(t)}^{-1}X-
E(Q_eTR_{W^{(\sigma)}(t)}^{-1}X)\] in $T_e{\cal D}^s_\mu({\cal
T}^n)$). Evidently, for $X=u(t)$, we obtain \[\breve
U_{u(t)}(t,m)=u(t,m-\sigma w(t))-U(t,m)\] and so
\[u(t,m-\sigma w(t))=U(t,m)+\breve U_{u(t)}(t,m)\] and
$E\breve U_{u(t)}(t,m)=0$. Then one can easily see that
\[E([(u\cdot\nabla)u](t,m-\sigma w(t)))=(U\cdot\nabla)U+E[(\breve
U_{u(t)}\cdot\nabla)\breve U_{u(t)}].\] Thus, \eqref{Reynolds}
transforms into
\begin{equation}\label{Reynolds2}
{\frac{\partial}{\partial
t}}U+(U\cdot\nabla)U-{\frac{\sigma^2}{2}}\nabla^2 U-{\rm
grad}\,p=-E[(\breve U_{u(t)}\cdot\nabla)\breve U_{u(t)}]
\end{equation}
which is the standard form of Reynolds equation. It differs from the
Navier-Stokes type relation with viscosity ${\frac{\sigma^2}{2}}$ by
the external force $-E[(\breve U_{u(t)}\cdot\nabla)\breve U_{u(t)}]$
that depends on $u(t,m)$, not on $U(t,m)$. Recall that
\eqref{Reynolds2} describes turbulent motion of fluid if the
dependence of $E[(\breve U_{u(t)}\cdot\nabla)\breve U_{u(t)}]$ on
$U$ is given (say, derived from experimental data, see
\cite{Sedov}).

\begin{remark}\label{4.1}
Note that for $\xi_t(s)$ introduced above, formula \eqref{D*=U}
reads that $\bar D_*\xi_t(s)_{|s=t}=U(T-t)$. Then, taking into
account formula \eqref{(8.12)}, one can easily derive that \[\bar
PD_*D_*\xi_t(s)_{s=t}=\bar
PD_*U(T-s,\xi_t(s))_{s=t}={\frac{\partial}{\partial
t}}U+(U\cdot\nabla)U-{\frac{\sigma^2}{2}}\nabla^2 U-{\rm grad}\,p.\]
Thus, eq. \eqref{Reynolds2} implies that $\bar
PD_*D_*\xi_t(s)_{s=t}=-PE[(\breve U_{u(T-t)}\cdot\nabla)\breve
U_{u(T-t)}]$.
\end{remark}

Our next aim is to show that a slight modification of the above
scheme of arguments allows us to annihilate the external force in
\eqref{Reynolds2} by introducing  a special random force field on
${\cal D}^s_\mu({\cal T}^n)$ into  \eqref{Euler1}.

For a random divergence free a.s. $H^{s+1}$-vector field
$X_\omega(m)$ on ${\cal T}^n$ (i.e., for a random vector
$X_\omega\in T_e{\cal D}^{s+1}_\mu({\cal T}^n)\subset T_e{\cal
D}^{s}_\mu({\cal T}^n)$), construct the random vector field $\breve
U_{X_\omega}(t,m)$ which, for any $\omega\in\Omega$, is given by the
formula
\[\breve U_{X_\omega}(t,m)=X_\omega(m-\sigma
w_\omega(t))-E(X_\omega(m-\sigma w_\omega(t)).\] Introduce the
non-random $H^s$ vector field $PE[(\breve
U_{X_\omega}\cdot\nabla)\breve U_{X_\omega}]$ and then construct the
random vector ${\frak F}_\omega(t,X_\omega)$ in $T_e{\cal
D}^s_\mu({\cal T}^n)$ by the formula
\[{\frak
F}_\omega(t,X_\omega)=Q_eTR_{W^{(\sigma)}_\omega(t)}PE[(\breve
U_{X_\omega}\cdot\nabla)\breve U_{X_\omega}].\] Note that
$PE[(\breve U_{X_\omega}\cdot\nabla)\breve U_{X_\omega}]$ and so
${\frak F}_\omega(t,X_\omega)$ lose the derivatives, i.e., they are
$H^s$-vector fields only since $X_\omega$ (and so $\breve
U_{X_\omega}$) is $H^{s+1}$. Thus ${\frak F}_\omega(t,X_\omega)$ is
well-posed only on an everywhere dense subset $T_e{\cal
D}^{s+1}_\mu({\cal T}^n)$ in $T_e{\cal D}^s_\mu({\cal T}^n)$.

Now introduce the right-invariant force vector field $\bar{\frak
F}_\omega(t,g,Y_\omega)$, where $Y_\omega\in T_g{\cal D}^s_\mu({\cal
T}^n)$, on ${\cal D}^s_\mu({\cal T}^n)$ that at $g\in {\cal
D}^{s+1}_\mu({\cal T}^n)$ and $\omega\in\Omega$ is determined by the
formula \[\bar{\frak F}_\omega(t,g,Y_\omega)=TR_g{\frak
F}_\omega(t,TR_g^{-1}Y_\omega),\] where $TR_g^{-1}Y_\omega$ is a
divergence free a.s. $H^{s+1}$-vector field.

Consider the equation
\begin{equation}\label{Newtonaux}
{\frac{\widetilde D}{dt}}\dot g_\omega(t)=\bar{\frak
F}_\omega(t,g_\omega(t),\dot g_\omega(t))
\end{equation}
on ${\cal D}^s_\mu({\cal T}^n)$ whose right hand side is well-posed
on the everywhere dense subset ${\cal D}^{s+1}_\mu({\cal T}^n)$ in
${\cal D}^s_\mu({\cal T}^n)$. Note that \eqref{Newtonaux} has no
diffusion term and so it is an ordinary differential equation with
parameter $\omega\in\Omega$. Here we do not investigate solvability
of \eqref{Newtonaux} but suppose that for the initial condition
$g_\omega(0)=e$ and $\dot g_\omega(0)=u_0\in T_e{\cal
D}^{s+1}_\mu({\cal T}^n)$, it a.s. has a unique $H^{s+1}$-solution
$g_\omega(t)$ which is a.s. well-posed on a non-random time interval
$t\in[0,T]$ for a certain $T>0$. Consider the divergence free a.s.
$H^{s+1}$-vector field $u_\omega(t,m)$ on ${\cal T}^n$ given by the
relation $\dot g_\omega(t)= u_\omega(t,g_\omega(t))$. The analog of
above-mentioned vector $U$ now takes the form
\begin{equation}\label{Uform}
{{\mathbb U}}(t,m)= E\bigl(u_\omega(t,m-\sigma w_\omega(t))\bigr)=
EQ_eTR_{W^{(\sigma)}_\omega(t)}^{-1}u_\omega(t).
\end{equation}
As well as in Lemma \ref{Udivfree} it is easy to see that vector
field \eqref{Uform} is divergence free.

\begin{theorem}\label{main1}
The divergence free vector field ${{\mathbb U}}$ given by
\eqref{Uform}, satisfies the Navier-Stokes equation without external
force and with viscosity $\frac{\sigma^2}{2}$:
\begin{equation}\label{NS}
{\frac{\partial}{\partial t}}{{\mathbb U}}+({{\mathbb
U}}\cdot\nabla){{\mathbb U}}-{\frac{\sigma^2}{2}}\nabla^2 {{\mathbb
U}}-{\rm grad}\,p_1=0.
\end{equation}
\end{theorem}

{\bf Proof.} Note that for the random field $u'_\omega(t,m)$ of
linear operators and the random field $u''_\omega(t,m)$ of bilinear
operators (the primes denote derivatives of $u$ in $m\in {\cal
T}^n$) the stochastic integrals $\int_0^t
u'_\omega(t,m)dw_\omega(t)$ and $\int_0^t
u''_\omega(t,m)(dw_\omega(t),dw_\omega(t))=\int_0^t {\rm
tr}\,u_\omega''dt=\int_0^t \nabla^2 u_\omega\, dt$ are well-posed.
Then by applying standard arguments to the Taylor series expansion
of $u_\omega$, one can easily see that the It\^o formula is
well-posed for $u_\omega(t,m-\sigma w_\omega(t))$ and so
\begin{gather*}
E\bigl(du_\omega(t,m-\sigma w_\omega(t))\bigr)=
E\left({\frac{\partial}{\partial t}}u_\omega (t,m-\sigma
w_\omega(t))dt+{\frac{\sigma^2}{2}}\nabla^2 u_\omega(t,m-\sigma
w_\omega(t))dt\right). \end{gather*} From \eqref{Newtonaux} it
follows (see \eqref{Euler1} and \eqref{Euler2}) that
${\frac{\partial}{\partial
t}}u_\omega=-P[(u_\omega\cdot\nabla)u_\omega]+{\frak
F}_\omega(t,u_\omega(t))$. Thus, in the same manner as in the proof
of Theorem \ref{mainaux} and deriving \eqref{Reynolds2}, we obtain
\begin{gather}
{\frac{\partial}{\partial t}}{{\mathbb U}}(t,m)=E({\frac{d}{d
t}}u_\omega(t,m-\sigma w_\omega(t)))=\notag\\
-E\Bigl[\Bigl((u_\omega\cdot\nabla)u_\omega\Bigr)(t,m-\sigma
w_\omega(t))\Bigr]
+{\frac{\sigma^2}{2}}\nabla^2 {{\mathbb U}}+{\rm grad}\,p\notag\\
+ EQ_eTR_{W^{(\sigma)}(t)}^{-1}{\frak
F}_\omega(t,u_\omega(t))=\label{iiii}\\ -({{\mathbb
U}}\cdot\nabla){{\mathbb U}} +{\frac{\sigma^2}{2}}\nabla^2 {{\mathbb
U}}+{\rm grad }p-E[(\breve {U}_{u_\omega(t)}\cdot\nabla)\breve
{U}_{u_\omega(t)}]\notag\\+ EQ_eTR_{W^{(\sigma)}(t)}^{-1}{\frak
F}_\omega(t,u_\omega(t))\notag.
\end{gather}
But by construction and by formulae \eqref{**i} and \eqref{**ii} we
get
\begin{gather}
EQ_eTR_{W^{(\sigma)}_\omega(t)}^{-1}{\frak F}_\omega(t,u_\omega(t))=\notag\\
EQ_eTR_{W^{(\sigma)}_\omega(t)}^{-1}Q_eTR_{W^{(\sigma)}_\omega(t)}PE[(\breve
{U}_{u_\omega(t)}\cdot\nabla)\breve {U}_{u_\omega(t)}]=\notag\\
EQ_eTR_{W^{(\sigma)}_\omega(t)}^{-1}TR_{W^{(\sigma)}_\omega(t)}Q_{W^{(\sigma)}_\omega(t)^{-1}}PE[(\breve
{U}_{u_\omega(t)}\cdot\nabla)\breve {U}_{u_\omega(t)}]=\notag\\
PE[(\breve {U}_{u_\omega(t)}\cdot\nabla)\breve
{U}_{u_\omega(t)}]\label{443}.
\end{gather}

Recall that for any divergence free fields ${\mathbb U}$, the vector
fields ${\frac{\partial}{\partial t}}{{\mathbb U}}$ and
$\nabla^2{{\mathbb U}}$ are divergence free. Hence, ${\rm grad}\,p$
in (\ref{iiii}) is taken from relation \eqref{P1} for
${E([(u_\omega\cdot\nabla)u_\omega](t,m-\sigma w_\omega(t)))}$,
i.e.,
\[PE([(u_\omega\cdot\nabla)u_\omega](t,m-\sigma
w_\omega(t)))=E([(u_\omega\cdot\nabla)u_\omega](t,m-\sigma
w_\omega(t)))-{\rm grad}\,p.\] Introduce ${\rm grad}\,p_1$ and ${\rm
grad}\,p_2$ by relations
\begin{gather*}
P({{\mathbb U}}\cdot\nabla){{\mathbb U}}=({{\mathbb
U}}\cdot\nabla){{\mathbb U}}-{\rm grad}\,p_1\\
\hspace{-155mm}\text{~~ and}\\
PE[(\breve {U}_{u_\omega(t)}\cdot\nabla)\breve
{U}_{u_\omega(t)}]=E[(\breve {U}_{u_\omega(t)}\cdot\nabla)\breve
{U}_{u_\omega(t)}]-{\rm grad}\,p_2.
\end{gather*}
Evidently, ${\rm
grad}\,p={\rm grad}\,p_1+{\rm grad}\,p_2$ (i.e., to within additive
constants $p=p_1+p_2$). Thus \eqref{NS} follows form \eqref{iiii}
and \eqref{443} in the natural form ${\frac{\partial}{\partial
t}}{{\mathbb U}}+({{\mathbb U}}\cdot\nabla){{\mathbb
U}}-{\frac{\sigma^2}{2}}\nabla^2 {{\mathbb U}}-{\rm grad}\,p_1=0.$
$\Box$


\end{document}